\input amstex
\input pstricks

\documentstyle{amsppt}
\newcount\headno
\newcount\subheadno 
\newcount\formno
\newcount\cno

\def\nextheadno{\global\advance\headno by 1 \the\headno}
\def\nextformno{\global\advance\formno by 1 \the\formno}
\def\nextno{\global\advance\cno by 1 \the\cno}

\def\headnum{\nextheadno. }
\def\procnum{\nextno}
\def\eqnum{\tag{\nextformno}}

\def\proclabel#1{\edef #1{\the\cno}}
\def\procref#1{$#1$}
\def\eqlabel#1{\edef #1{\the\formno}}
\def\eqref#1{($#1$)}
\def\headlabel#1{\edef #1{\the\headno}}

\def\CarD{|\partial\Delta|}
\def\Z{\Bbb Z}
\def\ZN{$\Bbb Z^d$}

\def\E{\Bbb E}
\def\R{\Bbb R}

\def\N{\Bbb N}
\def\LL{\Cal L}
\def\W{\Cal W}
\def\F{\Cal F}
\def\U{\Cal U}
\def\SD{\Bbb S}

\def\cl{\text{cl}}
\def\pc{\widehat{p_c}}
\def\scr{\scriptstyle}
\def\noti{\noindent\hskip 6mm}
\def\pxyc{\widetilde p_{xy}^{\,c}}
\def\pxyd{\widetilde p_{xy}^{\,d}}
\def\tvi{\vrule height 14pt depth 5pt width 0pt}
\def\setd{\mathbin \triangle}
\def\FBx{\Cal F_\Lambda^{B_{\scr x}}}

\newcount\refno
\global\refno=0
\def\nextrefno{\global\advance\refno by 1 }
\nextrefno\edef\Aba{\the\refno}
\nextrefno\edef\AlexI{\the\refno}
\nextrefno\edef\AlexII{\the\refno}
\nextrefno\edef\AlexIII{\the\refno}
\nextrefno\edef\AlexIV{\the\refno}

\nextrefno\edef\ACC{\the\refno}
\nextrefno\edef\Arr{\the\refno}
\nextrefno\edef\Bodineau{\the\refno}
\nextrefno\edef\Cerf{\the\refno}
\nextrefno\edef\CerfPiszI{\the\refno}
\nextrefno\edef\CerfPiszII{\the\refno}
\nextrefno\edef\CerfII{\the\refno}
\nextrefno\edef\Cha{\the\refno}
\nextrefno\edef\CouII{\the\refno}
\nextrefno\edef\Ferr{\the\refno}
\nextrefno\edef\GrimmII{\the\refno}
\nextrefno\edef\Grimm{\the\refno}
\nextrefno\edef\Grim{\the\refno}
\nextrefno\edef\Maes{\the\refno}
\nextrefno\edef\Rem{\the\refno}
\nextrefno\edef\Lig{\the\refno}
\nextrefno\edef\Pisz{\the\refno}

\topmatter
\NoRunningHeads
\title Poisson approximation for large clusters in the supercritical
FK
model
\endtitle
\author Olivier COURONN\'E\footnote{\hbox{Universit\'e Paris X-Nanterre, 
 Equipe modal'X,
        92001 Nanterre Cedex, France}\hfill\break
e--mail\;:\;couronne\@clipper.ens.fr\hfill\hss}
\endauthor
\date
13 May 2005
\enddate
\thanks
I thank R. Cerf for suggesting the problem and for critically reading
the manuscript. 
\endthanks

\abstract
Using the Chen-Stein method, we show that the spatial distribution
of large finite clusters in the supercritical FK model 
approximates a Poisson process when the ratio weak mixing property holds. 

\endabstract

\endtopmatter
\document
\subhead{\headnum Introduction}\endsubhead

We consider here the behaviour of
large finite clusters in the supercritical FK model.
In dimension two and more, 
their typical structure 
is described by the Wulff shape
\cite{\AlexIII, \AlexIV, \ACC, \Bodineau, \Cerf, \CerfPiszI, \CerfPiszII}. 
An interesting issue is the spatial distribution
of these large finite clusters.
Because of their rarity, a Poisson process
naturally comes to mind.
Indeed, we prove that the point process of the mass centers of
large finite clusters sharply approximates a Poisson process.
Furthermore, considering
large finite clusters in a large box such that their mean number is 
not too large, 
we observe Wulff droplets
distributed according to
this Poisson process. 

Redig and Hostad have recently studied the law of large finite clusters in 
a given box \cite{\Rem}. Their aim was different, in that they
obtained
accurate estimates on the law of the maximal 
cluster in the box, but the intermediate steps are similar. 
In the supercritical regime they considered 
only 
Bernoulli percolation and not FK percolation.

As in \cite{\Aba, \Cha, \Ferr, \Rem}, our main result is based on a 
second moment inequality. We have to 
control the interaction between two clusters. For this, 
we suppose that ratio weak mixing holds \cite{\AlexI}.
This property allows us to apply the Chen-Stein method in order to get
the approximation by a Poisson process. 

The ratio weak 
mixing holds in dimension two as soon as dual connectivities are
exponentially decreasing \cite{\AlexI}. For dimensions at least
three, we prove that ratio weak mixing holds for $p$ close enough to $1$.
Hence our main results are valid in all dimensions for $p$ large enough.

The following section 
is devoted to the statement of our results.
In section 3, we define the FK model. We recall the weak and the ratio 
weak mixing properties and we state a perturbative mixing result in section 4.
Section 5 contains the definition of our point process and 
the description of the Chen-Stein method.
The core of the article is section 6, where we study a second moment
inequality.
In section 7, we deal with the probability of having a large finite cluster
with its center at the origin.
In section 8, we treat the case of distant clusters and we finish the proof
of Theorem 1. 
The proof of Theorem 3 is done in section 9, and the proof of the
perturbative mixing result is done in section 10.

\subhead{\headnum Statement of the results}\endsubhead
We consider the FK measure $\Phi$ on the $d$-dimensio\-nal lattice \ZN\ 
and in the supercritical regime. The point $\pc$ stands for $\widehat
p_g$ in 
dimension two, and for $p_c^{\text{slab}}$ in dimensions three and
more. For $q\geq 1$ we let $\Cal U(q)$ be the set such that there
exists a unique FK measure on \ZN\ of parameters $p$ and $q$ if $p$ is not in
$\Cal U(q)$. By \cite{\Grimm} this set is at most countable.

Let $\Lambda$ be
a large box in \ZN. We fix $n$ an integer and we consider the finite clusters
of cardinality larger than $n$.
We call them $n$-{\it large clusters\/}.
Let $C$ be a finite cluster. The {\it mass center} of $C$ is
 $$M_C=\bigg\lfloor\frac{1}{|C|}\sum_{x\in C}x\bigg\rfloor,$$
where $\lfloor x\rfloor$ denotes the
site of \ZN\  whose coordinates are the integer part of those of $x$.
We define a 
process $X$ on $\Lambda$ by
$$X(x)=\bigg\{\eqalign{&1 \hbox{ if }x\hbox{ is the mass center of a $n$--large cluster }C\cr    
 &0 \hbox{ otherwise.}}$$

Let $\lambda$ be the expected number of sites $x$ in $\Lambda$
such that $X(x)=1$. We denote by $\LL(X)$ the law of a process
$X$. For $Y$ a process on $\Lambda$, 
we let $||\LL(X)-\LL(Y)||_{TV}$ be the total variation distance between
the laws of the processes $X$ and $Y$ \cite{\Arr}.

\proclaim{Theorem \procnum}
Let $q\geq 1$ and $p>\pc$ with $p\notin\Cal U(q)$.
Let $\Phi$ be the FK measure on \ZN\ of parameters $p$ and $q$.
We suppose that $\Phi$ is ratio weak mixing. There exists a constant
$c>0$ such 
that:
for any box $\Lambda$, letting $X$ be defined as
above, and letting $Y$ 
be a Bernoulli process on $\Lambda$ with the same one-dimensional marginals 
as $X$,
we have for $n$ large enough
$$||\LL(X)-\LL(Y)||_{TV}\leq \lambda\exp(-cn^{(d-1)/d}).$$
\endproclaim
\proclabel{\prin}

As a corollary, the number of large clusters in $\Lambda$
is approximated by a Poisson variable.
\proclaim{Corollary \procnum}
Let $\Phi$ be as in Theorem~\procref{\prin}. Let $N$ be the number of
large finite clusters 
whose mass centers are in the box $\Lambda$. Let $Z$ be a Poisson variable of
mean $\lambda$, and let $c>0$ be the same constant as in
Theorem~\procref{\prin}.  
Then for any $A\subset \Z^+$ and for $n$ large enough,
 $$|P(N\in A)-P(Z\in A)|\leq\lambda\exp\big(-cn^{(d-1)/d}\big).$$
\endproclaim
\proclabel{\coro}

We provide next a control of the shape of the large finite
clusters. Here we consider a sequence of boxes $(\Lambda_n)_n$. If the size of 
$\Lambda_n$ is not too large, that is of order less than $\exp(\rho n^{(d-1)/d})$ 
for a certain constant $\rho$, then the energy created by 
the $n$--large clusters 
of $\Lambda_n$ dominates a term of entropy. In this case we can assert that 
the shape of these $n$-large clusters are close to the Wulff shape.

More precisely, let $\W$ be the Wulff crystal, let $\theta$ be the 
density of the 
infinite cluster, and let $\LL^d(\cdot)$ be the Lebesgue measure on
$\R^d$. Let
$$W=\frac{1}{\big(\theta \LL^d(\W)\big)^{1/d}}\W$$
be the renormalized Wulff crystal. For $l>0$, let
$V_\infty(C,l)$ be the  
neighbourhood of $C$  of width $l$ for the metric $|\cdot|_\infty$.
For two sets $A$ and $B$, the notation $A\setd B$ stands for the symmetric
difference between $A$ and $B$.

\proclaim{Theorem \procnum}
Let $\Phi$ be as in Theorem~\procref{\prin}. Let $f:\N\rightarrow\N$
be such
that $f(n)/n\rightarrow 0$ and $f(n)/\ln n\rightarrow\infty$ 
as $n$ goes to infinity. Let $(\Lambda_n)_n$ be a sequence of boxes in \ZN, 
and let $\lambda_n$ be the expected number of mass centers
of $n$--large clusters in 
$\Lambda_n$. 
For all $\delta>0$, there exists $c>0$ such that if
$\limsup 1/n^{(d-1)/d}\ln \lambda_n\leq c$.
$$\eqalign{\limsup_{n\rightarrow\infty}\frac{1}{n^{(d-1)/d}}\ln\Phi\Big[
\LL^d\Big(&\big(\bigcup_{\scr x\in \Lambda_n\atop \scr X(x)=1}
(x+W) 
\big) \setd\cr
& \big(n^{-1}\bigcup_{\scr C\ n\text{-large}\atop \scr
M_C \in\Lambda_n} 
V_\infty(C, f(n))\big)\Big)
\geq \delta\big|\{x:X(x)=1\}\big|\Big]<0.}$$
\endproclaim
\proclabel{\deuze}
\noindent
For clarity, we omit the subscript $n$ on $X$.

\noindent
{\bf Remark :} Consider a sequence $(\Lambda_n)_n$ such that 
$|\Lambda_n|\simeq \exp(\rho n^{(d-1)/d})$ and let $w_1>0$ be such 
that~\cite{\CerfII}:
$$P(n\leq|C(0)| <\infty)\approx \exp(-w_1n^{(d-1)/d}).$$
On the one hand we need 
$\rho\geq w_1$ in order to have some $n$--large clusters in $\Lambda_n$. 
On the other hand the condition on $\lambda_n$ 
in theorem~\procref{\deuze} may be 
rewritten as $\rho\leq c+w_1$.

\smallskip
The ratio weak mixing property is a key hypothesis in our results. 
The following proposition allow us to apply the three preceding results for
$p$ large enough in all dimensions.

\proclaim{Proposition \procnum}
Let $d\geq3$ and $q\geq1$. There exists $p_0<1$ such that $\Phi$ satisfies the ratio
weak mixing property for $p>p_0$.
\endproclaim
\proclabel{\tres}

\subhead{\headnum FK model }\endsubhead

We consider the lattice \ZN\  with $d\geq2$. We turn it into a graph
by adding bonds between all pairs $x$, $y$   
of nearest neighbours. We write $\E$ for the set of bonds and we let
$\Omega$ be the set $\{0,1\}^\E$. A 
{\it bond configuration\/} $\omega$ is 
an element of 
$\Omega$. A bond $e$ is {\it open\/} in $\omega$ if $\omega(e)=1$,
and {\it closed\/}  
otherwise. 

A {\it path\/} is a sequence $(x_0, \dots, x_n)$ of
distinct sites 
such that $\langle x_i, x_{i+1}\rangle$
is a bond for each $i$, $0\leq i\leq n-1$. 
A subset $\Delta$ of \ZN\ is {\it connected\/} if for every $x$, $y$
in $\Delta$, there exists a path included in $\Delta$ connecting $x$
and $y$.
If all bonds of a path are open in 
$\omega$, we say that the path is {\it open} in $\omega$.
A {\it cluster\/} 
is a connected component in \ZN\ when we keep only open bonds. It is
usually denoted by $C$.  
Let $x$
be a site. We write $C(x)$ for  
the cluster containing $x$.

To define the FK measure, we first consider finite volume FK measures.
Let $\Lambda$ be a box included in \ZN. 
We write $\E(\Lambda)$ for the set of bonds 
$\langle x,y\rangle$
with $x,y\in\Lambda$. 
Let $\Omega_\Lambda=\{0,1
\}^{\E(\Lambda)}$ be the space 
of bonds configuration in $\Lambda$. 
Let $\Cal F_\Lambda$ be its $\sigma$-field, that is the set of subsets
of $\Omega_\Lambda$.
For $\omega$ in $\Omega_\Lambda$, we define $\cl(\omega)$ as the
number of clusters 
of the configuration $\omega$. 

For $p\in[0,1]$ and $q\geq 1$,
the FK measure in $\Lambda$ with parameters $p,q$ and free boundary 
condition is the probability measure on $\Omega_\Lambda$ 
defined by
$$\forall\, \omega \in \Omega_\Lambda \quad
\Phi_\Lambda^{f,p,q}(\omega)=\frac{1}{Z_\Lambda^{f,p,q}}
\Big(\prod_{e\in\E(\Lambda)}p^{\omega (e)}(1-p)^{1-\omega (e)}\Big)
q^{\cl(\omega)},$$
where $Z_\Lambda^{f,p,q}$ is the appropriate normalization factor.

We also define FK measures for arbitrary boundary conditions.
For this, let $\partial\Lambda$ be the {\it boundary} of $\Lambda$,
$$\partial\Lambda=\{x\in\Lambda \text{ such that } \exists\,y\notin\Lambda,
\langle x,y\rangle \text{ is a bond}\}.$$
For a partition $\pi$ of $\partial\Lambda$, a $\pi$--cluster is a
cluster of $\Lambda$ when we add open bonds between the pairs of sites
that are in the same class of $\pi$. 
Let
$\cl_\pi(\omega)$ be the number of $\pi$--clusters in $\omega$.
To define $\Phi^{\pi,p,q}_\Lambda$ we replace
$\cl(\omega)$ by $\cl_\pi(\omega)$ and $\smash{Z_\Lambda^{f,p,q}}$ by
$\smash{Z_\Lambda^{\pi,p,q}}$ in the above formula. 

There exists
a countable subset $\U(q)$ in $[0,1]$ such that the following holds.
As $\Lambda$ grows
and invades the whole lattice \ZN, the finite volume measures converge
weakly toward the same infinite measure $\Phi_\infty^{p,q}$ 
for all $p\notin \U(q)$ \cite{\Grimm}. 
We will always suppose that this occurs, that is 
$p\notin \U(q)$.
We shall drop the superscript and the 
subscript on $\Phi_\infty^{p,q}$, and simply write $\Phi$.  
It is known that the FK measure $\Phi$ is translation--invariant.

The measure $\Phi$ verify the {\it finite energy property\/}: for each $p$
in $(0,1)$, 
there exists $\delta>0$ such that for every finite--dimensional cylinders 
$\omega_1$
and $\omega_2$ that differ by only one bond, 
$$\Phi(\omega_1)/\Phi(\omega_2)\geq\delta.\eqnum$$
\eqlabel{\finite} 

The random cluster model has a phase transition. There exists
$p_c\in(0,1)$ such that there is no infinite cluster $\Phi$--almost
surely if $p<p_c$, and an infinite cluster $\Phi$--almost
surely if $p>p_c$. Other critical points have been introduced
in order to work with 'fine' properties.
In dimension two, we define $\widehat p_g$ as the critical point for
the exponential  
decay of dual connectivities, see \cite{\CouII, \Grimm}. 
In three and more dimensions, let $p_c^{\text{slab}}$ be 
the limit of the critical points for the percolation in slabs \cite{\Pisz}.
For brevity, $\widehat p_c$ will stand for $\widehat p_g$ in dimension two,
and for $p_c^{\text{slab}}$ in dimensions three and more. It is
believed that $\pc=p_c$ in all dimensions and for all $q\geq 1$, 
but in most cases we know only that $\pc\geq p_c$.

We now state Theorem $17$ of \cite{\CerfII}, applied to FK measures.

\noindent
If $q\geq 1$, $p>\widehat p_c$ and $p\notin\Cal U(q)$, there exists 
$w_1>0$ such that
$$\lim \frac{1}{n^{(d-1)/d}}\ln\Phi\big(n\leq|C(0)|<\infty\big)=-w_1,\eqnum$$
\eqlabel{\limit}

\noindent
where $C(0)$ is the cluster of the origin.

\subhead{\headnum Mixing properties}\endsubhead

Let $x$ and $y$ be two points in \ZN\ and let $(x_i)_{i=1}^d$ and 
$(y_i)_{i=1}^d$ be their coordinates.  
Write $|x-y|_1=\smash{\sum_{i=1}^d}|x_i-y_i|.$ 

\proclaim{Definition \procnum}
Following \cite{\AlexII},
we say that $\Phi$ satisfies the \rom{weak mixing property} 
if for some $c, \mu > 0$,
for all sets $\Lambda, \Delta \subset \Z^d$,
$$\eqalign{
  \sup \big\{\big|&\Phi(E \mid F) - \Phi(E)\big|:  E \in \F_{\Lambda}, F \in
    \F_{\Delta}, \Phi(F) > 0\big\} \cr
  &\leq c \sum_{x \in \Lambda,y \in \Delta} 
    e^{-\mu |x - y|_{\scr 1}}. }\eqnum$$
\eqlabel{\weakmix}
\endproclaim

\proclaim{Definition \procnum}
Following \cite{\AlexII},
we say that $\Phi$ satisfies the \rom{ratio weak mixing property} 
for some $c_1, \mu_1 > 0$,
for all sets $\Lambda, \Delta \subset \Z^d$,
$$\eqalign{ 
  \sup &\Big\{ \Big| \frac{\Phi(E \cap F)}{\Phi(E)\Phi(F)} - 1 \Big| : E \in 
    \F_{\Lambda}, F \in
    \F_{\Delta}, \Phi(E)\Phi(F) > 0 \Big\} \cr
  &\leq c_1 \sum_{x \in \Lambda,y \in \Delta}
    e^{\scr -\mu_{\scr 1} |x - y|_{\scr 1}}, }\eqnum$$
\eqlabel{\ratiomix}
\endproclaim

\noindent
Roughly speaking, the influence of what happens in $\Delta$ on 
the state of the bonds
in $\Lambda$ decreases exponentially with the distance between $\Lambda$
and $\Delta$. 

In dimension two, the measure $\Phi$ is ratio
weak mixing as soon as 
$p>\widehat p_c$ \cite{\AlexII}, but such a result is not available in
dimension larger than three. 
We provide a perturbative mixing result, which is valid for all
dimensions larger than three, and which is similar to the weak mixing
property.

\proclaim{Lemma \procnum}
Let $d\geq 3$ and $q\geq 1$. There exists $p_1<1$ and $c>0$ such that:
for all $p>p_1$, all connected sets $\Gamma,\Delta$
with $\Gamma\subset\Delta$,
every boundary conditions $\eta$, $\xi$  
on $\Delta$,
every event $E$ supported on $\Gamma$,
$$|\Phi_\Delta^{\eta,p,q}(E)-\Phi_\Delta^{\xi,p,q}(E)|\leq 
2\CarD
\exp\big(-c\inf\big\{|x-y|_1, x\in\Gamma, y\in\partial\Delta\big\}\big)
.$$
\endproclaim
\proclabel{\disagree} 
\noindent
We are not aware of a particular reference for this result,
and we give a sketch of the proof in Section 10.

\subhead{\headnum The Chen-Stein method}\endsubhead

From the percolation process, we want to extract a point process
describing the occurrence of large finite clusters.
For a point $x$ in $\R^d$, let $\lfloor x\rfloor$ denotes the
site of \ZN\  whose coordinates are the integer parts of those of $x$.
Assume that $C$ is a finite subset of \ZN. Then the {\it mass center}
of $C$ is 
$$M_C=\bigg\lfloor\frac{1}{|C|}\sum_{x\in C}x\bigg\rfloor.$$
Let $n\in\N$.
A $n$--{\it large cluster} is a finite cluster of cardinality larger
than $n$. Let $\Lambda$ be a box in \ZN.
We define a 
process $X$ on $\Lambda$ by
$$X(x)=\bigg\{\eqalign{&1 \hbox{ if }x\hbox{ is the mass center 
of a $n$--large cluster }C\cr   
 &0 \hbox{ otherwise.}}$$

In order to apply the Chen-Stein method, we define for $x, y$ in \ZN,
$$\eqalign{p_x=&\ \Phi(X(x)=1),\cr
p_{xy}=&\ \Phi\big(\exists\, C,C'\text{ two clusters such that: }
C\cap C'=\emptyset,\cr
&\quad\quad \ \ n\leq |C|,|C'|<\infty, M_C=x\text{ and }M_{C'}=y
\big),}$$ and we let
$B_x=B(x, n^2)$ be the box centered at $x$ of side length $n^2$.
Let $\lambda$ be the expected number of sites $x$ in $\Lambda$ 
such that $X(x)=1$. We have $\lambda=\sum_{x\in \Lambda}p_x$ and,
because of the translation--invariance of $\Phi$, for each site $x$ in
$\Lambda$ 
$$\lambda=|\Lambda|\cdot p_x.\eqnum$$
\eqlabel{\first}

We introduce three coefficients $b_1$, $b_2$, $b_3$ by:
$$b_1=\sum_{x\in\Lambda}\sum_{y\in B_x}p_xp_y,$$
$$b_2=\sum_{x\in\Lambda}\sum_{y\in B_x\setminus x}p_{xy},$$
$$b_3=\sum_{x\in\Lambda}E\Big|E\Big(X(x)-p_x|\sigma(X(y), y\notin
B_x \big)\Big)\Big|.$$

\noindent
Let $Z_1$ and $Z_2$ be two Bernoulli processes on $\Lambda$. 
The {\it total variation distance}
between the laws of the processes $Z_1$ and $Z_2$ \cite{\Arr} is
$$||\LL(Z_1)-\LL(Z_2)||_{TV}=2\sup\big\{\big|P(Z_1\in A)-P(Z_2\in A)\big|,A
\text{ subset  
of } \{0,1\}^\Lambda\big\}.$$

\noindent
Let $Y$ be a Bernoulli process on $\Lambda$ such that the $Y(x)$'s are
iid and $$P(Y(x)=1)=p_x.$$
The Chen-Stein method provides a control of the total variation
distance between $X$ and $Y$ in terms of the $b_i$'s. Indeed we apply  
Theorem~$3$ of \cite{\Arr} to obtain that
$$||\LL(X)-\LL(Y)||_{TV}\leq 2(2b_1+2b_2+b_3)+4\sum_{x\in\Lambda} p_x^2.
\eqnum$$
\eqlabel{\chen}

To prove Theorem~\procref{\prin}, we shall provide an upper bound
on each term $b_i$.
The ratio weak mixing property is essential to our proof of 
the bound of
$b_2$. Nevertheless, we believe that one can prove the
following inequality, without any mixing assumption:
$$\Phi\big[n\leq C(x)<\infty, n\leq C(y)<\infty, C(x)\cap C(y)=\emptyset\big]
\leq\Phi(2n\leq C(0)<\infty).\eqnum$$

Let us give now an upper bound on $p_x$.
By \cite{\GrimmII}, there exists a constant $c>0$ such that:
$$\Phi(n\leq|C(0)|<\infty)\leq\exp\big(-cn^{(d-1)/d}\big).$$
But
$$\eqalign {p_x\leq&\sum_{k\geq n}\Phi\big(\exists\, C, |C|=k, M_C=x\big)\cr
\leq&\sum_{k\geq n}\sum_{y\in B(x, 2k)}\Phi\big(|C(y)|=k\big)\cr
\leq&\sum_{k\geq n}(2k)^d\exp\big(-ck^{(d-1)/d}\big).}$$
Hence there exists a constant $c>0$ such that for $n$ large enough
$$p_x\leq\exp(-cn^{(d-1)/d}).\eqnum$$
\eqlabel{\px}

\noindent

\subhead{\headnum Second moment inequality}\endsubhead
In this section we bound the term $p_{xy}$ with the help of the
ratio weak mixing property. First we introduce a local version
of $p_{xy}$.
We define $\widetilde
p_{xy}$ by
$$\eqalign{\widetilde p_{xy}=\Phi\big(&\exists\, C, C'\text{ two clusters} 
\text{ such that }\cr
&n\leq|C|<n^2, n\leq|C'|<n^2,
M_{C}=x, \text{ and }M_{C'}=y\big).}$$

The distance between
two sets $\Gamma$ and $\Delta\subset\Z^d$ is
$$d(\Gamma,\Delta)=\inf\{|x-y|_1,x\text{ in }\Gamma, y\text{ in }\Delta\},$$
and it is the length of the shortest path connecting $\Gamma$ to $\Delta$.

We divide the term $\widetilde p_{xy}$ into two parts. 
Let $\mu_1$ be the constant appearing in the definition of the ratio
weak mixing property and let $K>5/\mu_1$. We define 
$\pxyc$ by 
$$\eqalign{\pxyc=\Phi\big(&\exists\, C, C'\text{ two clusters} 
\text{ such that }d(C, C')\leq K\ln n,\cr
&n\leq|C|<n^2, n\leq|C'|<n^2,
M_{C}=x, \text{ and }M_{C'}=y\big).}$$
We define also $\pxyd$ by
$$\eqalign{\widetilde p_{xy}=\Phi\big(&\exists\, C, C'\text{ two clusters} 
\text{ such that }d(C, C')> K\ln n,\cr
&n\leq|C|<n^2, n\leq|C'|<n^2,
M_{C}=x, \text{ and }M_{C'}=y\big).}$$
The superscripts $c$ and $d$ stand for {\it close} and {\it distant}.
So $\widetilde p_{xy}=\pxyc+\pxyd$
and we study separately these two terms.

First we focus on $\pxyd$. We have
$$\eqalign{
\pxyd\leq
\sum_{C, C' \text{ distant}}\Phi(C 
\text{ and }C'\text{ are clusters}),}$$
where the sum is over the couples $(C, C')$ of connected subsets of \ZN\ such
that $$\eqalign{n\leq|C|<n^2, &\ n\leq|C'|<n^2,\cr 
&M_C=x, M_{C'}=y, \text{ and }d(C, C')>K\ln n.}$$
Let $c_1,\mu_1$ be the constants appearing in the definition of the ratio
weak mixing property. Let $(C,C')$ be a couple appearing
in the sum above. We have
$$\sum_{u \in C,v \in C'}
    e^{-\mu_1 |u - v|}\leq n^4\exp(-\mu_1K\ln n),$$
so for $n$ large enough 
$$c_1 \sum_{u \in C,v \in C'}
    e^{-\mu_1 |u - v|}\leq 1.$$ 
So for $n$ large enough
$$\Phi(C \text{ and }C'\text{ are clusters})\leq 
2\Phi(C\text{ is a cluster})\cdot\Phi(C'\text{ is a cluster}),$$
by the ratio weak mixing property \eqref{\ratiomix}. 
Hence there exists $c>0$ such that for $n$ large enough

$$\eqalign{\pxyd&\leq\sum_{u\in B(x, 2n^2), v\in B(y, 2n^2)}
2\Phi(n\leq|C(u)|<\infty)\cdot\Phi(n\leq|C(v)|<\infty)\cr
&\leq \exp(-cn^{(d-1)/d}).}\eqnum$$
\eqlabel{\pxydl}

Now we consider $p_{xy}^{\,c}$. We have

$$\eqalign{
\pxyc\leq
\sum_{C, C' \text{ close}}\Phi(C \text{ and }C'\text{ are clusters}),}$$
where the sum is over the couples $(C, C')$ of subsets of \ZN\ such
that 
$$\eqalign{n\leq|C|<n^2, &n\leq|C'|<n^2, \cr
&M_C=x, M_{C'}=y, \text{ and }d(C, C')\leq K\ln n.}$$ 
For $n$ large enough, the event $\{C \text{ and }C'\text{ are clusters}\}$ 
is $\Cal F_{B(x,3n^2)}$-measurable. 
So we only consider bond configurations in
$B(x,3n^2)$. 

We give a deterministic total order on the pairs $(u,v)$ of \ZN\ in
such a way that  
if $|u_1-v_1|_1<|u_2- v_2|_1$, then $(u_1, v_1)<(u_2, v_2)$.
Let $(C,C')$ be a pair of sets appearing in the above sum.
Take a configuration 
$\omega$ in $B(x,3n^2)$ such that $C$ and $C'$ are clusters in $\omega$.
We change the configuration
$\omega$ as follows. 

To start with, we take the pair $(u,v)$ such that
$u\in C$, $v\in C'$ and $(u,v)$ is the first such pair for the order 
above. For $0\leq i\leq d$, we define $t_i$ the point whose $d-i$ first
coordinates are equal to those of $u$, and the others are equal to
those of $v$. Hence $t_0=u$, $t_d=v$, and $t_i$ and $t_{i+1}$ differ
by only one coordinate.
We consider the shortest path $(u_0,\dots, u_k)$ connecting $u$ to $v$
through the $t_i$'s. It is composed of the segments $[t_i, t_{i+1}]$
for $0\leq i\leq d-1$.

We open all the bonds $\langle u_i, u_{i+1}\rangle$ for 
$i=0\dots k-1$. In the same time, we close all the bonds
incident to $u_i$ for $i=1\dots k-1$ distinct from the previous 
bonds $\langle u_j, u_{j+1}\rangle$. 
Let $\widetilde \omega$ be the new configuration
in $B(x,3n^2)$. 
We denote by $\widetilde C$ the set $C\cup C'\cup \{u_i\}_{i=1}^{k-1}$.
By construction,  $\widetilde C$ is a cluster in $\widetilde \omega$.
We have
$$2n\leq \widetilde C<4n+K\ln n.$$
The number of bonds we have changed is bounded by 
$2dK\ln n$. By the finite energy property \eqref{\finite}:
$$\Phi(\widetilde\omega)\geq n^{2dK\ln\delta}\Phi(\omega),$$
for a certain constant $\delta$ in $(0, 1)$.

Now we control the number of antecedents by our transformation.
Take a configuration $\widetilde\omega$ of $B(x,3n^2)$.
To get an antecedent of $\widetilde\omega$, we have to

\noti(a) choose two sites $u,v$ in $B(x,3n^2)$, with $|u-v|_1\leq K\ln
n$

\noti(b) take the path connecting $u$ to $v$ along the coordinate axis

\noti(c) choose the state of the bonds that have an endpoint on this path.

In step $(a)$, we have less than $(3n^2)^d(2K\ln n)^d$ choices. 
In step $(b)$ we have
just one choice. In step $(c)$ the number of choices is bounded by 
$2^{2dK\ln n}$.  
Hence for $n$ large enough 
the number of antecedents of $\widetilde\omega$ is bounded 
by $n^{4dK}$. 

Finally, 
$$\sum_{C, C' \text{ close}}\Phi(C 
\text{ and }C'\text{ are clusters})
\leq n^{4dK}\cdot n^{2dK\ln\delta}\sum_{\widetilde C}\Phi(
\widetilde C \text{ is a cluster}),$$
where the sum is over connected subsets $\widetilde C$ of \ZN\ such that
$2n\leq |\widetilde C|<5n$ and $\widetilde C$ is contained in
$B(x,3n^2)$. This sum is bounded by 
$$|B(x,3n^2)|\cdot\Phi(2n\leq |C(0)|<5n).$$
Thus by \eqref{\limit}, there exists $c_2>w_1$ such that for $n$ large enough,
$$\widetilde p_{xy}^{\,c}\leq \exp(-c_2n^{(d-1)/d}).\eqnum$$
\eqlabel{\pxycl}

To conclude, remark that
$$p_{xy}-\widetilde p_{xy}\leq \Phi\big(\exists\, C \text{ a cluster
such that }  
n^2\leq|C|<\infty, M_C=x\big).$$
By \eqref{\px}, 
there exists $c$ such that for $n$ large enough 
the difference 
between $p_{xy}$ and $\widetilde p_{xy}$ is bounded
by $\exp(-cn^{2(d-1)/d})$.
So by \eqref{\pxydl} there exists $c>0$ such that
$p_{xy}\leq\pxyc+\exp(-cn)$. Since in \eqref{\pxycl} the constant
$c_2$ is strictly larger than $w_1$, 
there exists $c_3>w_1$ such that
for $n$ large enough 
$$p_{xy}\leq \exp(-c_3n^{(d-1)/d}).\eqnum$$ 

\eqlabel{\moment}

\subhead{\headnum A control of $p_x$}\endsubhead
We compare $p_x$ and $\Phi(n\leq|C(0)|<\infty)$.
\proclaim{Lemma \procnum}
If $q\geq1$, $p>\widehat p_c$, and
$p\notin\Cal U(q)$, then 
$$\lim \frac{1}{n^{(d-1)/d}}\ln p_x=-w_1.$$
\endproclaim
\proclabel{\ratiolim}

\noindent
We note that in \cite{\Rem}, the authors take the left endpoints of 
the clusters instead of the mass center and get the same limit.

\demo{Proof of Lemma~\procref{\ratiolim}} We begin with a lower bound
for $p_x$. We recall that for all $x$ in \ZN, 
$p_x=\Phi(X(0)=1)$. Let $\alpha>1$. Because  of \eqref{\limit}, we have
$$\lim \frac{1}{n^{(d-1)/d}}\ln\Phi(n\leq |C(0)|<\infty)= \lim
\frac{1}{n^{(d-1)/d}}\ln \Phi(n\leq |C(0)|<n^\alpha).$$ Then 
$$\eqalign{
\Phi(n\leq |C(0)|<n^\alpha)&\leq \sum_{x\in B(0, n^\alpha)}
\Phi(n\leq|C(0)|<n^\alpha, M_C=x) \cr 
&\leq |B(0, n^\alpha)|\Phi(X(0)=1).
}$$

We give next an upper bound:
$$\eqalign{\Phi(X(0)=1)&=\Phi(\exists C \text{ a cluster}, M_C=0,
n\leq|C|<n^\alpha) \cr
&\quad\quad\quad\quad\quad{}+\Phi(\exists C \text{ a cluster}, M_C=0,
n^\alpha\leq|C|<\infty)\cr
&\leq\smash{\sum_{x\in B(0, n^\alpha)}}\tvi\Phi(n\leq|C(x)|<\infty) \cr
&\quad\quad\quad\quad\quad+ \sum_{k\geq n^\alpha}\Phi\big(\exists C \text{ a
cluster}, |C|=k, C\cap 
B(0, 2k)\neq\emptyset\big)\cr
&\leq |B(0, n^\alpha)|\Phi(n\leq|C(x)|<\infty)+
\sum_{k\geq n^\alpha} |B(0, 2k)|\Phi(|C(0)|=k).
}$$

Finally, we use the limit \eqref{\limit} to get 
$$\lim\frac{1}{n^{(d-1)/d}}\ln p_x=\lim \frac{1}{n^{(d-1)/d}}
\ln\Phi(n\leq|C(0)|<\infty)=-w_1.\qed$$
\enddemo

\subhead{\headnum Proof of Theorem~\procref{\prin}}\endsubhead

We recall that $\Lambda$ is a box and $\lambda$ is the expected number
of the mass centers in $\Lambda$ of $n$--large clusters.
We write $\F_\Lambda^{B_{\scr x}}$ for the $\sigma$--field
$\F_{\Lambda\setminus B_{\scr x}}$.
First, we bound the term
$$E\big|E\big(X(x)-p_x|\FBx\big)\big|.$$

\noindent
Let $\widetilde X(x)$ be equal to $1$ if $x$ is the mass center of a cluster
$C$, with $C$ such that $n\leq|C|< n^2/4$, and equal to $0$ otherwise.
Let $\widetilde p_x=\Phi(\widetilde X(x))$. We have
$$\eqalign{E\big|E\big(X(x)-p_x|\FBx&\big)\big|
\leq E\big|E\big(X(x)-\widetilde X(x)|
\FBx\big)\big|\cr
{}+&E\big|E\big(\widetilde X(x)-\widetilde p_x|
\FBx\big)\big|
+E\big|E\big(\widetilde p_x-p_x|\FBx\big)\big|.}
\eqnum$$
\eqlabel{\trois} 

\noindent
Since the quantity $X(x)-\widetilde X(x)$
is always positive, 
$$\eqalign{E\big|E\big(X(x)-\widetilde X(x)|
\FBx\big)\big|&=E\big [E\big(X(x)-\widetilde X(x)|
\FBx\big)\big ]\cr
&=p_x-\widetilde p_x.}$$
We have also 
$$E\big|E\big(\widetilde p_x-p_x|\FBx\big)\big|=p_x-\widetilde p_x.$$
But
$$p_x-\widetilde p_x=\Phi(\exists\, C\text{ a cluster}, n^2/4\leq |C|<\infty,
M_C=x),$$
so by \eqref{\px} there exists $c>0$ such that
$p_x-\widetilde p_x \leq\exp(-cn^2).$

The variable $\widetilde X(x)$
is $\F_{B(x,n^2/4)}$-measurable. 
The distance between $B(x, n^2/4)$
and the complementary region of $B_x$ is of order $n^2$.
If $\Phi$ is weak mixing,
or by lemma~\procref{\disagree} if $p$ is close enough to $1$, there
exists a constant $c>0$ such 
that for $n$ large enough 
$$E\big|E\big(\widetilde X(x)-\widetilde p_x|
\F_\Lambda^{B_x}\big)\big|\leq \exp(-cn^2).$$

Putting together the estimates of the three terms on the 
right-hand side of \eqref{\trois}, we conclude that there
exists $c>0$ such that for $n$ large enough
$$E\big|E\big(X(x)-p_x|\F_\Lambda^{B_x}\big)\big|\leq \exp(-cn^2).\eqnum$$
\eqlabel{\btrois}

Now observe that
$|\Lambda|=\lambda p_x^{-1}.$ 
Using inequality \eqref{\moment} and the limit of
Lemma~\procref{\ratiolim}, there exists $c>0$ such that 
$$b_2\leq \lambda p_x^{-1} \exp\big(-c_3n^{(d-1)/d}\big)\leq \lambda 
\exp\big(-cn^{(d-1)/d}\big).$$

Because of \eqref{\btrois}, there exists $c>0$, $c'>0$ such that
$$b_3\leq \lambda p_x^{-1} \exp(-cn^2)\leq  \lambda\exp(-c'n^2).$$

The term $b_1$ is controlled by Lemma~\procref{\ratiolim}.
We apply finally the Chen-Stein inequality 
\eqref{\chen} to obtain
Theorem~\procref{\prin}.\qed

\subhead{\headnum Proof of Theorem~\procref{\deuze}}
\endsubhead

The Wulff crystal is the typical shape of a large finite cluster
in the supercritical regime. The crystal is built on
a surface tension $\tau$. 
The surface tension is a function from $\SD^{d-1}$,
the $(d-1)$--dimensional unit sphere of $\R^d$, to $\R^+$.
It controls the exponential decay of the probability
for having a large separating surface in a certain direction, 
with all bonds closed.
We refer the reader to \cite{\Cerf, \CerfII}
for an extended
survey of this function. 

In the regime
$p>\widehat p_c$ and $p\notin\Cal U(q)$, the surface tension is positive,
continuous, and satisfies the weak simplex inequality.
We denote by $\W$ the Wulff shape associated to $\tau$,
$$\W=\{x\in \R^d, x.u\leq\tau(u) \text{ for all }u\text{ in }\SD^{d-1}\}.$$
The Wulff shape is a main ingredient in the proof of \eqref{\limit}. 

Let $\theta=\Phi(0\leftrightarrow\infty)$ be the density of the infinite
cluster. Let $f:\N\rightarrow\N$, such
that $f(n)/n\rightarrow 0$ and $f(n)/\ln n \rightarrow\infty$ 
as $n$ goes to infinity. 
Let $x$ and $y$ be two points of $\R^d$, and let $(x_i)_{i=1}^d$ and 
$(y_i)_{i=1}^d$
be their coordinates. We write $|x-y|_\infty=\max_{1\leq i\leq d} |x_i-y_i|$.
We define a neighbourhood of a cluster $C$ by
$$V_\infty(C,f(n))=\{x\in \R^d,\exists\, y \in C,|x-y|_\infty\leq f(n)\}.$$

Let $(\Lambda_n)_n$ be a sequence of boxes in \ZN, 
and let $\lambda_n$ be the expected number of mass centers
of $n$--large clusters in 
$\Lambda_n$.
In Theorem~$3$, we consider the event
$$\eqalign{\bigg\{\LL^d\Big(&\bigcup_{\scr x\in\Lambda_n \atop \scr X(x)=1}
(x+\theta\LL^d(\W)^{-1/d}\W 
\big) \setd\cr
&n^{-1}\bigcup_{\scr C \ n\text{--large}\atop \scr
M_C \in\Lambda_n} 
V_\infty(C, f(n))\Big)
\geq \delta\big|\{x:X(x)=1\}\big|\bigg\}.}\eqnum$$
\eqlabel{\gros}
It is included in the event
$$\eqalign{\bigg\{&\text{there exists }C \text{ a }n\text{--large cluster}
\text{ such that } M_C\in\Lambda_n,\cr
&\LL^d\Big(\big(
M_C+\theta\LL^d(\W)^{-1/d}\W 
\big) \setd
\big(n^{-1}
V_\infty(C, f(n))\big)\Big)
\geq \delta\bigg\}.}$$
Taking the logarithm of its probability and dividing by $n^{(d-1)/d}$,
we may show that for $n$ large it is equivalent to the logarithm divided by
$n^{(d-1)/d}$ of the following
quantity: 
$$\lambda_n\Phi\Big[\LL^d\Big(\big(
M_{C(0)}+\theta\LL^d(\W)^{-1/d}\W 
\big) \setd
\big(n^{-1}
V_\infty(C(0), f(n))\big)\Big)
\geq \delta\big| n\leq|C(0)|<\infty\Big].$$
By \cite{\Cerf, \CerfII}, there exists $c>0$ such that if
$$\limsup
1/n^{(d-1)/d}\ln \lambda_n\leq c,$$ 
then the inequality in Theorem~\procref{\deuze} holds. \qed

\subhead{\headnum A perturbative mixing result}\endsubhead

First we prove lemma~\procref{\disagree}, following the proof of the
uniqueness of the FK measure 
for $p$ close enough to $1$ in \cite{\Grim}. The difference is that we
consider not just one but   
two independent FK measures. The idea of using two 
independent copies of a measure comes from \cite{\Maes}. Then the
proof of proposition~\procref{\tres} follows.

\demo{Proof of lemma~\procref{\disagree}}

Let $\Delta$ be a connected subset of \ZN. There is a partial order
$\preceq$ in $\Omega_\Delta$ given by 
 $\omega\preceq\omega^\prime$ if and only if $\omega(e)\leq\omega^{\prime}(e)$
 for every bond $e$. A function $f:\Omega_\Delta\rightarrow\R$ is called
 {\it increasing} if $f(\omega)\leq f(\omega^{\prime})$ whenever
 $\omega\preceq\omega^{'}$. An event is an element of
$\Omega_\Delta$. An event is called increasing if its 
 characteristic function is increasing. 
For a pair of probability
 measures $\mu$ and $\nu$ on $(\Omega_\Delta, \Cal F_\Delta)$, 
we say that $\mu$
 {\it(stochastically) dominates} $\nu$ if for any $\Cal F_\Delta$-measurable
 increasing function $f$ the expectations satisfy
 $\mu(f)\geq\nu(f)$ and we denote it by $\mu\succeq\nu$.
Let $P_p$ be the Bernoulli bond--percolation measure on \ZN\ of parameter $p$. 
The FK measures on $\Delta$ dominate stochastically a certain
Bernoulli measure 
restricted on $\E(\Delta)$:
$$\Phi_\Delta^{\eta,p,q}\succeq
P_{p/[p+q(1-p)]}\big|_{\E(\Delta)}. \eqnum$$  
\eqlabel{\bern}

For $(\omega_1,\omega_2)\in \Omega^2$, we call a site $x$ {\it white}
if $\omega_1(e)\omega_2(e)=1$ for all bond $e$ incident with $x$,
and {\it black} otherwise.
We define a new graph structure on \ZN.
Take two sites $x$ and $y$ and label $x_i, y_i$ their coordinates. 
If $\max_{i=1 ... d}|x_i-y_i|=1$, 
then $\langle x,y \rangle$ is a {\it $\star$-bond} and $y$ is
a {\it $\star$-neighbour} of $x$. 
A {\it $\star$-path} is a sequence $(x_0, ..., x_n)$ of distinct sites
such that 
$\langle x_i, x_{i+1}\rangle$ is a $\star$-bond for $0\leq i\leq
n-1$.
 
For any set $V$ of sites, the {\it black cluster} $B(V)$ is the union
of $V$ together with the set of all $x_0$ for which there exists a
$\star$-path $x_0,\dots, x_n$ such that $x_n\in V$ and 
$x_0,\dots, x_{n-1}$ are all black.
Let $\Gamma$, $\Delta$ be two connected sets
with $\Gamma\subset\Delta$. 
The 
'interior boundary' $D(B(\partial\Delta))$ of $B(\partial\Delta)$ is the set
of sites $x$ 
satisfying:
\goodbreak

\noti(a)\ $x\notin B(\partial\Delta)$

\noti(b)\ there is a $\star$-neighbour of $x$ in $B(\partial\Delta)$

\noti(c)\ there exists a path from $x$ to $\Gamma$ that does not use 
a site in $B(\partial\Delta)$.

\noindent
Let $I$ be the set of sites $x_0$ 
for which there exists a path $x_0,\dots,x_n$ with $x_n\in\Gamma$,
$x_i\notin B(\partial\Delta)$ for all $i$, see figure $1$.
\midinsert
\vbox{\centerline{\hbox{\pspicture(-5,-4.2)(5,4.3)
\psframe[fillcolor=cyan,fillstyle=solid,linewidth=2pt](-3.6,-3.6)(3.6,3.6)
\pscurve[fillcolor=white,fillstyle=solid](-2.9,-0.1)(-2.6,-2.8)(0,-3)(2,-2.3)(3,-1.7)(2.4,0,2)(2,1.6)(2.3,2.7)
(1.2,2.6)(-0.2,2.4)(-1.8,2.6)(-2.6,2.4)(-2.5,1.2)(-2.9,-0.1)
\psframe[fillcolor=white,fillstyle=solid,linewidth=2pt](-1.1,-1.1)(1.1,1.1)
\rput(-0.2,-0.1){$\Gamma$}
\rput(-4.9,1.48){$\Delta$}
\psline[linewidth=0.5pt]{->}(-4.7,1.4)(-3.6,1.2)
\rput(-5.5,-0.4){$D(B(\partial\Delta))$}
\psline[linewidth=0.5pt]{->}(-4.6,-0.4)(-2.9,-0.1)
\rput(4.7,-0.9){$I$}
\psline[linewidth=0.5pt]{->}(4.5,-0.9)(2.2,-1)
\endpspicture
}}
\centerline{figure $1$: The set $I$ inside $\Delta$}
}
\endinsert
\noindent
Let $$K_{\Gamma, \Delta}= ·\big\{\big(B(\partial\Delta)\cup
D(B(\partial\Delta))\big)  
\cap\Gamma=\emptyset\big\}.$$  
If $K_{\Gamma, \Delta}$ occurs, we have the following facts:

\noti(a)\ $D(B(\partial\Delta))$ is connected

\noti(b)\ every site in $D(B(\partial\Delta))$ is white

\noti(c)\ $D(B(\partial\Delta))$ is measurable with respect to the colours 
of sites in $\Z^d\setminus I$

\noti(d)\ each site in $\partial I$ is adjacent to some site 
of $D(B(\partial\Delta))$.

\noindent
These claims have been established in the proof of 
Theorem $5.3$ in \cite{\Grim}.
 
Pick $\eta,\xi$ two boundary conditions of $\Delta$.
For brevity let $\Cal P=\Phi_\Delta^{\eta,p,q}\times \Phi_\Delta^{\xi,p,q}$.
We shall write $X, Y$ for the two projections from
$\Omega_\Delta\times\Omega_\Delta$ to $\Omega_\Delta$. 
Then for any $E\in\Cal F_\Gamma$, we have by the claims above 
$$\Cal P(X\in E,K_{\Gamma, \Delta})=\Cal P(Y\in E,K_{\Gamma, \Delta})=
\Cal P(\Phi^{w,p,q}_I(E)1_{K_{\Gamma,\Delta}}).$$
Hence $$|\Phi_\Delta^{\eta,p,q}(E)-\Phi_\Delta^{\xi,p,q}(E)|\leq \big(1-\Cal
P(K_{\Gamma, \Delta})\big).$$ 

Because of inequality \eqref{\bern} and by the stochastic
domination result in \cite{\Lig}, 
the process
of black sites is stochastically dominated by a Bernoulli
site--percolation 
process whose parameter is independent of $\Gamma$,
$\Delta$, $\eta$, $\xi$ and decreases to $0$ as $p$ goes to $1$. 
There exists $p_1<1$ such that this Bernoulli process is 
subcritical for the $\star$-graph structure of \ZN\ and for $p\geq p_1$. 
Hence there exists $c>0$ such that for $p>p_1$, for all $\Gamma$, $\Delta$,
$\eta$, $\xi$,

$$\Cal P(K_{\Gamma, \Delta})\geq 1-\CarD\exp\big(-c\,d(\Gamma, 
\partial \Delta)\big).\quad\qed$$ 
\enddemo

\demo{Proof of proposition~\procref{\tres}}
The domination inequality
\eqref{\bern} implies that for $p$ large enough, the mesures 
$\Phi_\Delta^{\eta,p,q}$ have exponentially bounded controlling
regions in the terminology of~\cite{\AlexI}. Thus by theorem~$3.3$ of
\cite{\AlexI}, 
lemma~\procref{\disagree} implies the ratio weak mixing property for
the mesures $\Phi_\Delta^{\eta,p,q}$. \qed  
\enddemo

\end